\documentclass[11pt,a4paper, reqno]{amsart}
\usepackage{amsmath}
\usepackage[latin9]{inputenc}
\usepackage{amsfonts} 
\usepackage{amssymb}
\usepackage{color}
\allowdisplaybreaks
\usepackage{enumerate}
\usepackage[]{hyperref}

\newtheorem{thm}{Theorem}[section]
\newtheorem{cor}[thm]{Corollary}

\newtheorem{pro}[thm]{Proposition}
\newtheorem{dfn}[thm]{Definition}
\newtheorem{rem}[thm]{Remark}

\numberwithin{equation}{section}

\def\N{\mathbb N}
\def\Z{\mathbb Z}
\def\R{\mathbb R}

\theoremstyle{plain}
\title[Young Inequalities]{New Young inequalities and applications}

\author[Fernández-Martínez]{Pedro Fernández-Martínez}
\address[Pedro Fernández-Martínez]{Departamento de Matemáticas\\
Facultad de Matemáticas\\ Universidad de Murcia \\ Campus de Espinardo\\
30071 Espinardo (Murcia), Spain}
\email{pedrofdz@um.es}
\author[da Silva]{Eduardo Brandani da Silva}
\address[Eduardo Brandani da Silva]{Universidade Estadual de Maringá - UEM \\ Departamento de Matemática\\
Av. Colombo 5790\\ Maringá, Paraná \\ 870300-110, Brazil }
\email{ebsilva@wnet.com.br}
\thanks{This work was supported by the Spanish Ministerio de Economía y Competitividad (MTM2013-42220-P) and 
Fundaci\'on S\'eneca de la Regi\'on de Murcia 19378/PI/14.}

\subjclass{Primary 46B70, 47B07; Secondary 46E30}

\begin{document}
\maketitle
\begin{abstract}
We establish upper bounds for the convolution operator acting between interpolation spaces. This will provide several examples of Young Inequalities in different families of function spaces. We use this result to prove a bilinear interpolation theorem and we show applications to the study of bilinear multipliers.
\end{abstract}
\section{Introduction}
The real interpolation method introduced by Lions and Peetre in 1964, see \cite{LionsPeetre}, has proved to be a very useful tool in many areas of analysis such as harmonic analysis, partial differential equations, approximation theory, operator theory or functional analysis. See the monographs by Butzer and Berens \cite{Butzer1967}, Bergh and Löfström \cite{Bergh-Lofstrom}, Triebel \cite{Triebel1978,Triebel1992,Triebel2006}, Beauzamy \cite{Beauzamy1978}, König \cite{Koenig1986}, Bennett Sharpley \cite{Bennett-Sharpley}, Tartar \cite{Tartar2007} or the monographs by Connes \cite{Connes1994}, and Amrein, Boutet de Monvel and Geourgescu \cite{Amrein1996} for applications to other areas.

However the classical version of this method fails to  identify the end-point spaces of the interpolation scales it generates. As an example, let us recall that the classical real method does not produce Lorentz-Zygmund spaces from the couple $(L_1, L_{\infty})$.  In order to do so, we need to introduce limiting interpolation methods such as 
logarithmic  methods.  The papers by Evans and Opic \cite{EO} and also  Evans, Opic and Pick \cite{EOP}, where the authors study logarithmic interpolation methods, inspired the appearance of the some other limiting methods defined by means of slowly varying functions and rearrangement invariant (r.i.) spaces. These have been studied by T. Signes and one the the present authors in  \cite{FMSS,FMS-1,FMS-2,FMS-4,FMS-3,FMS-5} and allow to produce limit spaces that are not in the classical real interpolation scale.

On the other hand, several papers that study bilinear interpolation theorems have been recently published. See for example Mastylo \cite{Mast2013} or Cobos and Segurado \cite{CSe-Log} where we can find bilinear interpolation theorems for logarithmic methods. Here, we extend the study of bilinear  interpolation theorems to the methods  defined by slowly varying functions and r.i. spaces. 
The main obstacle for our approach is the lack of a Young type inequality for r.i. spaces. For this reason,  in a first stage we  establish a general Young inequality in the context of r.i. spaces that will be used to prove subsequent bilinear interpolation theorems. 
In order to be more precise, let the measure space $(\Omega, \mu)$  be $\R$ with the Lebesgue measure or $\Z$ with the counting measure.
Given $f$ and $g$, measurable functions on  $(\Omega, \mu)$, 
the convolution of $f$ and $g$, $f \ast g$, is defined as 
\[ (f \ast g) (x) = \int_{\Omega} f(x-y) g(y) d\mu . \]
The classical 
Young inequality estimates the norm of the convolution  $ f \ast g $ in $ L^{r} $ when  $ f  \in  L^{p} $ and $ g  \in  L^{q} $ with
$ \frac{1}{p} + \frac{1}{q} =1 + \frac{1}{r}   $ . In fact, 
 \[ \| f \ast g \|_{L^{r}}  \leq \| f \|_{L^{p}} \| g \|_{L^{q}} .\]
This inequality was extended by O'Neil to the context of Lorentz spaces in \cite{ONeil1963} and to Orlicz spaces in \cite{ONeil1965}. Further generalizations of Young inequality were carried out in the following years. See the papers by Hunt \cite{Hunt1966}, Yap \cite{Yap1969}, and Blozinski  \cite{Blozinski72AMS}, \cite{Blozinski1972}.
More recent contributions studying Young Inequality in the context of weighted Lebesgue spaces are due to Kerman \cite{Kerman1983}, and Bui \cite{BUI1994}. The inequality in Lebesgue spaces with variable exponent was studied by Samko in \cite{Samko1998-2} and
\cite{Samko1998-1}. 

The study of the boundedness of the convolution operator uses a variety of different techniques. Here we obtain estimates for the norm of the convolution operator acting among  r.i. spaces by using  a simple interpolation argument that already appears in \cite{Bergh-Lofstrom} for $L_p$ spaces.
  The  inmediate applications of this result recover most of  the known estimates for the convolution operator. For example, if we consider the convolution operator acting between Lebesgue spaces we recover the classical Young inequality, but we also obtain Young inequalities  for Lorentz spaces or  the inequalities for Orlicz spaces established by O'Neil in \cite{ONeil1965}. We can also recover some of the recently established inequalities for weighted Lorentz spaces, see the papers by Krepela in \cite{Kvrepela2014}, \cite{Krepela2014} 
 and \cite{Kvrepela2016}.

Once we have made a convolution inequality available, we are in a position to tackle the problem of establishing bilinear interpolation theorems for the interpolation methods defined by means of slowly varying functions and r.i. spaces. This is done in  \S \ref{section-bilinear}.
We close the paper with  \S \ref{bilinear multipliers} where we show some  applications of these results  to the  study Fourier multipliers. In the appropiate context, the use of our bilinear interpolation theorems improves some of the results published by Blasco in  \cite{Blasco2009}.

Throughout the paper, we use the notation $f \lesssim g$ to indicate that there exist a constant $\mathbf{c}>0$ such that the functions (or any other quantities depending on some parameters) $f$ and $g$ satisfy  $f \leq \mathbf{c} g$. We write $f \sim g$ if $f \lesssim g$ and $g \lesssim f$.

\section{Young Inequality}

Let  $E$ be a function space   over the measure space $(\Omega, \mu)$, let $f \in E$ while $g$ belongs to the associate  space $ E' $. The norm of the convolution $ f \ast g $ satisfies the inequalities
\begin{align}
& \|f \ast g \|_{E}  = \Big \| \int_{\Omega} f(x-y) g(y) d\mu \Big \|_{E} \leq \int_{\Omega} \|f\|_{E} |g(y)| d\mu \leq \|f\|_{E} \| g \|_{L^{1}} \label{conv I} \\
& \|f*g\|_{L^{\infty}} \leq \Big \| \int_{\Omega} f(x-y) g(y) d\mu \Big \|_{L^{\infty}} \leq \|f\|_{E} \| g \|_{E'}. \label{conv II}
\end{align} 
This enable us to  establish the following theorem.
\begin{thm} \label{GYI-thm}
Let $ \mathcal{F} $ be an exact interpolation functor and let $ E $ be a Banach function space, then for any $ f \in E $ and $ g \in E' $
\[ \| f \ast g \|_{\mathcal{F}(E, L^{\infty})} \leq \| f \|_{E} \; \|  g \|_{\mathcal{F}(L^{1}, E')} . \]  
\end{thm}
\begin{proof}
Fix $f \in E$ and consider the operator
\[ Tg = f \ast g . \]
Equations \eqref{conv I} and \eqref{conv II}
show that  the operators 
\begin{eqnarray*}
T: L^{1} & \longrightarrow & E \\
T: E' & \longrightarrow & L^{\infty}
\end{eqnarray*} 
are bounded with norms no greater than $ \|f\|_{E}$. Now, by using the interpolation functor $\mathcal{F}$ we obtain that
$$T: \mathcal{F}(L^{1}, E') \longrightarrow \mathcal{F}(E, L^{\infty})$$
with norm no greater than $ \|f\|_{E}$. In particular, for any $g \in E'$ 
\begin{equation} \label{GYI}
\|f \ast g\|_{\mathcal{F}(E, L^{\infty})} \leq  \| f \|_{E} \; \|g\|_{\mathcal{F}(L^{1}, E' )}. 
\end{equation}

\end{proof}
\begin{rem} Should the interpolation functor $\mathcal{F}$ not be exact then the convolution inequality \eqref{GYI} would be 
\[\|f \ast g\|_{\mathcal{F}(E, L^{\infty})} \leq  C \| f \|_{E} \; \|g\|_{\mathcal{F}(L^{1}, E' )}\]
where $C = \max\{ \| T \|_{L^{1}, E}, \| T \|_{E', L^{\infty}}  \}$.
\end{rem}
This simple interpolation argument, together with appropiate choices of the interpolation functor $ \mathcal{F} $ and the space $ E $, provides other forms of Young's inequality. Let us show some examples.

\begin{cor}
Let $L^{\varphi_{0}}$ be an Orlicz space on $\R$ with absolutely continuous norm, and let
$ L^{\varphi} $ be the Orlicz space whose Young function satisfies $ \varphi^{-1} =(\varphi_{0}^{-1})^{1 - \theta} .$ Then, for any measurable function $ g $
\[ \| f \ast g \|_{L^{\varphi}}  \leq \|f \|_{L^{\varphi_{0}}} \; \|  g \|_{L^{\psi}} ,\]
where 
$ \psi^{-1}(t) = t \Big( \dfrac{\psi_{0}^{-1}(t)}{t}  \Big)^{\theta} $ and $ L^{\psi} $ is the corresponding Orlicz space endowed with the Orlicz norm. Here $ \psi_{0} $ is the complementary Young function of $\varphi_{0}$.
\end{cor}
\begin{proof}
Consider the complex interpolation method $ [\cdot, \cdot]_{\theta} $, with $0 < \theta < 1$.
Theorem 1.14 of \cite{KPS} together with Theorems 7.3 and 9.1 of \cite{GP} yield, with equality of norms,  that
\begin{align*}
&[L^{\varphi_{0}}, L^{\infty}]_{\theta} = L^{\varphi} \\
&[L^{1}, L^{\psi_{0}}]_{\theta} = (L^{1})^{1 - \theta} (L^{\psi_{0}})^{\theta} = L^{\psi}.
\end{align*} 
Last equality holds with the same norm if we endowed $L^{\psi}$ with its Orlicz norm.
Now apply \eqref{GYI}
to obtain the result.
\end{proof}

We may obtain a more general inequality if we use Gustavsson-Peetre  method  to interpolate Orlicz spaces.
In order to do that we need to recall some concepts and facts from \cite{GP}.  
Given a  function $\varphi : (0, \infty)  \longrightarrow (0,\infty)$, the dilation function of $\varphi$ is defined as
 $s_{\varphi}(t) = \sup_{s >0} \frac{\varphi(ts)}{\varphi(s)}$ for $t>0$. 
The function $s_{\varphi}$ is submultiplicative, therefore we can define the lower and upper dilation indices of $\varphi$
as:
\begin{equation*}
\pi_{\varphi} = \sup_{t<1} \frac{\log s_{\varphi}(t)}{\log(t)} = \lim_{t \to 0}  \frac{\log s_{\varphi}(t)}{\log(t)},  \qquad 
\rho_{\varphi} = \sup_{t>1} \frac{\log s_{\varphi}(t)}{\log(t)} = \lim_{t \to \infty}  \frac{\log s_{\varphi}(t)}{\log(t)}.
\end{equation*}

Consider the Orlicz space $L^{\varphi_{0}}$, where the indices of  the function $\varphi_{0}$ satisfy
$0 < \pi_{\varphi_{0}} \leq \rho_{\varphi_{0}} < \infty$, and let $\rho:(0,\infty) \longrightarrow (0,\infty)$ be a pseudo-concave function. Then, if $s_{\rho}(t) = \circ (\max \{1, t\})$ for $t \to 0$ or $t \to \infty$ and we interpolate the couple $(L^{1}, L^{\varphi_{0}})$ using the Gustavsson-Peetre method with parameter $\rho$, we obtain
\begin{equation} \label{io1}
 \big < L^{1}, L^{\varphi_{0}}, \rho  \big > = L^{1} \, \rho \Big ( \frac{L^{\varphi_{0}}}{L^{1}} \Big ) =L^{\varphi} .
\end{equation}
Here $\varphi^{-1}(s) = s \, \rho \Big( \frac{\varphi_{0}^{-1}(s)}{s}    \Big)$, $s>0$. In case $s_{\rho} (t) \to 0$ as $t \to 0$
the Gustavsson-Peetre method applied to the couple $( L^{\varphi_{0}}, L^{\infty})$ produces the space
\begin{equation} \label{io2}
 \big <  L^{\varphi_{0}}, L^{\infty}, \rho  \big > = L^{\varphi_{0}} \rho \Big ( \frac{L^{\infty}}{L^{\varphi_{0}}} \Big )= L^{\varphi} .
\end{equation}

 Now, combining equations \eqref{GYI}, \eqref{io1} and \eqref{io2} we establish the following corollary.
 \begin{cor} Let $\varphi_{0}: (0,\infty) \longrightarrow (0, \infty)$ with dilation indices satisfying that 
 $0 < \pi_{\varphi_{0}} \leq \rho_{\varphi_{0}} < \infty$, and let $\rho:(0,\infty) \longrightarrow (0,\infty)$ be a pseudo-concave function with the property  $s_{\rho}(t) = \circ (\max \{1, t\})$ for $t >0$. Then
 \begin{equation*}
  \| f \ast g \|_{L^{\varphi_{0}} \rho \big( \frac{L^{\infty}}{L^{\varphi_{0}}}  \big)}  
 \leq \|f \|_{L^{\varphi_{0}}} \; \|  g \|_{  L^{1} \rho \big ( \frac{L^{\varphi_{0}}}{L^{1}} \big )    } .
 \end{equation*}
 \end{cor}

 We can also convolve periodic functions defined on the one-dimensional torus $\mathbb{T}$ with the usual measure. In this case the convolution  is defined as
 $$\big ( f \ast g \big)(e^{ix}) = \int_{0}^{2\pi} f(e^{i (x-y)}) g (e^{i y}) dy .$$ 
 Next result establishes a Young inequality for Lorentz-Zygmund spaces on the torus. 

\begin{cor}
 Let $f$ and $g$ be continuous functions on $\mathbb{T}$ with $f \in L \log L (\mathbb{T})$. Then, for $ 0 < \theta <1$, 
 $$ \|   f \ast g \|_{(L \log L)^{(1-\theta)} (L^{\infty})^{\theta} }  \lesssim \|   f \|_{L \log L}   \|   g \|_{(L^{1})^{1 - \theta} L_{exp}^{\theta}}. $$
\end{cor}
\begin{proof}
As before, the same arguments  apply here with the complex interpolation method. We also have to take into account that 
$(L \log L)' = L_{exp}$ with equivalence of norms, so we need to introduce $\lesssim$ in our inequality.
\end{proof}
 
 Our next result involves Lorentz-Karamata spaces. We begin with a definition that we take from \cite{FMS-2}. 
 Let $(\Omega,\mu)$ denote a  $\sigma$-finite measure space with a non-atomic measure $\mu$.
\begin{dfn}
Let $E=E(0,\infty)$ be an r.i. space, $1\leq p\leq \infty$ and let $b$ be a slowly varying function on $(0,\infty)$.
The  \textit{Lorentz-Karamata type} space $L_{p,b,E}=L_{p,b,E}(\Omega,\mu)$ consists of all $\mu$-measurable functions on $\Omega$ such that
\begin{equation}\label{L-K1}
\|f\|_{L_{p,b,E}}=\|t^{1/p}b(t)f^*(t)\|_{\widetilde{E}}<\infty.
\end{equation}
\end{dfn}
Here $\widetilde{E}$ stands for the r.i. space $E$ with respect to the homogeneous measure on $(0,\infty)$
$$\nu(A) = \int_{0}^{\infty} \chi_{A}(t) \frac{dt}{t}.$$
See \cite{FMS-1} \cite{FMS-2} and \cite{FMS-3} for a detailed description of $\widetilde{E}$.
When $p=\infty$, $L_{p,b,E}$ is different from the trivial space if, and only if, $\|b\|_{\widetilde{E}(0,1)}<\infty$. 

These spaces are particular examples of ultrasymmetric spaces studied by Pustylnik (see \cite{Pust-ultrasymm}). For the choice $E=L_q$, the space
$L_{p,b,E}$ coincides with the  Lorentz-Karamata spaces $L_{p,q;b}$. See \cite {EKP}, \cite{Neves2002} or \cite{GOT} for more information on Lorentz-Karamata spaces.

The family $L_{p,b,L_q}$ contains also well-known spaces for particular choices of $b$.
If $m\in\N$, $\alpha=(\alpha_1,\ldots,\alpha_m)$ and
$$b^{\alpha}(t)=\prod_{i=1}^{m} \ell_i^{\alpha_i}(t), \quad t>0$$
where
$$\ell_1(t)=\ell(t)=1+|\log t|,\quad \ell_i(t)=\ell(\ell_{i-1}(t)),\quad t>0,\quad i=2,\ldots,m,$$
then, the space $L_{p,b^\alpha,L_q}$ coincides with the \textit{generalized Lorentz-Zygmund} space $L_{p,q;\alpha}$ introduced by Edmunds, Gurka and Opic in \cite{ego1997} and studied in \cite{op}.
 Moreover, the space  $L_{p,b,L_q}$ when $b(t)=\ell^{\alpha}(t)$ is the \textit{Lorentz-Zygmund} space $L^{p,q}(\log L)^\alpha$ introduced by Bennett and Rudnick (see \cite{BR}, \cite{Bennett-Sharpley}). For $b\sim 1$, $L_{p,b,L_q}$
 coincides with Lorentz space $L^{p,q}$   and with the Lebesgue space $L_q$ if $b\sim 1$ and $p=q$.

Now we are in a position to state  the next three corollaries:
\begin{cor} \label{cor2.7}
 Let $f $ and $g$ be continuous functions on the torus $\mathbb{T}$ with $f$ in $L \log L (\mathbb{T})$. Then, for $0 < \theta < 1$ and $b$ any slowly varying function we have the estimate 
 $$  \|   f \ast g \|_{L_{\frac{1}{1 - \theta}, B_{\theta}, F}}    \lesssim \|   f \|_{L \log L}  \, \|   g \|_{L_{q, B_{\theta}, E}}, $$
 where $B_{\theta} (t) = \ell(t)^{- \theta} b (t \ell(t))$.
\end{cor}

\begin{cor} \label{cor2.8}
 For the case $\theta =0$ the following inequality holds, 
 $$  \|   f \ast g \|_{ L_{(1, B_{0}, F)}  \cap (L^{1}, L^{\infty})_{0, b(t \ell(t)), F, 1, L^{1}} }    \lesssim \|   f \|_{L \log L}  \, \|   g \|_{L_{(1,B_{0}, F)}} $$
\end{cor}

\begin{cor} \label{cor2.9}
 For the case $\theta =1$ the following inequality holds  
 $$  \|   f \ast g \|_{ L_{(\infty, B_{1}, F)  }}    \lesssim \|   f \|_{L \log L}  \, \|   g \|_{
L_{(1, B_{0}, F)}  \cap (L^{1}, L^{\infty})^{\mathcal{R}}_{1, b(t \ell(t)), F, \ell(t)^{-1}, L^{\infty}} } 
   $$
\end{cor}

\begin{proof}[Proof of Corollaries \ref{cor2.7}, \ref{cor2.8} and  \ref{cor2.9}]
In order to prove the above corollaries apply \eqref{GYI} and Corolaries 7.3 and 7.4  of \cite{FMS-2}.
\end{proof}

\section{Bilinear interpolation theorems} \label{section-bilinear}

In this section we use the estimate for the convolution operator \eqref{GYI} to establish  bilinear interpolation theorems for the interpolation methods defined by means of  slowly varying functions and rearrangement invariant spaces. We will work with $K$ and $J$-interpolation methods. 
%
Before we introduce them, let us recall the definitions of the $K$ and the $J$ functionals. 

We say $\overline{X}=(X_0,X_1)$ is a compatible couple of Banach spaces if  
$X_{0}$ and $X_{1}$ are Banach spaces  continuously embedded in the same
Hausdorff topological vector space. We equip $X_0+X_1$ with the norm $K(1,\cdot)$, where for $t>0$
$$K(t,f)=K(t,f;X_0,X_1)=\inf \big \{\|f_0\|_{X_0}+t\|f_1\|_{X_1}:\ f=f_0+f_1,\ f_i\in X_{i} \big \}$$
is the Peetre $K$-functional. We also consider the intersection  and the family of norms on  $X_0\cap X_1$ defined by the $J$ functional
\[ J(t, a; X_0,X_1)=\max\{\|f\|_{X_0},t\|f\|_{X_1}\}, \; t >0.\]
We refer to \cite{Bennett-Sharpley}, \cite{Bergh-Lofstrom} or \cite{B-K} for more information on interpolation theory.

The following definition describes the $K$-interpolation method defined with an slowly varying function and an r.i. espace. In particular, this extends the well known $K$-interpolation method 
with a function parameter introduced by Gustavsson, see~\cite{Gustfunctparam}. 

\begin{dfn} \label{dinter}
Let $\overline{X} = (X_{0}, X_{1})$ be a compatible Banach couple, $E$ a r.i.  space on $\R$,  $b$ a slowly varying function on $(0,\infty)$ and $0\leq\theta \leq 1$. The real interpolation space $\overline{X}^{K}_{\theta,b,E}
\sim (X_0,X_1)^{K}_{\theta, b, E}$ consists of all $f$ in $X_{0} + X_{1}$ for which the norm
$$ \|f\|^{K}_{\theta,b,E} = \big\| e^{-\theta t} {b}(e^{t}) K(e^{t},f)\big\|_{E}$$
is finite.
\end{dfn}

 This family of $K$-spaces comprises many other classes of interpolation spaces studied by different authors. Let us mention some of them.
 When $E=L_q$ and $b \sim 1$, the space $(X_0,X_1)_{\theta,b,E}$ coincides with the classical real interpolation spaces 
 $(X_0,X_1)_{\theta,q}$. Moreover,  since for $0 < \theta < 1$ the function $t^\theta/b(t)$ is equivalent to a function of the Kalugina class $B_{K}$, the spaces $(X_0,X_1)_{\theta,b,L_q}$  are special cases of the  interpolation spaces with a function parameter studied by Gustavsson in \cite{Gustfunctparam}. Interpolation spaces $(X_0,X_1)_{\theta,b,L_q}$, including limiting cases $\theta=0$ and $\theta=1$, have been studied by Gogatisvili, Opic and Trebels in \cite{GOT}, while $(X_0,X_1)_{\theta,\ell^{\mathbb A}(t),L_q}$ were thoroughly studied by Evans, Opic and Pick in \cite{EOP} and by Evans and Opic in \cite{EO}. For ordered couples, $X_0\hookrightarrow X_1$ and $\alpha\in\R$, the interpolation spaces $(X_0,X_1)_{\theta,(1+|\log t|)^\alpha,L_q}$ were defined and studied by Doktorskii in \cite{do}. More recently, and also for ordered couples,  Cobos, Fern\'andez-Cabrera, Kh\"un and Ullrich have studied the spaces $(X_0,X_1)_{0,\frac1{\log t},L_q}$  and have identified them with a new class of $J$-interpolation spaces denoted by $(X_0,X_1)_{0,q;J}$ when $1<q\leq\infty$ (see Theorem 4.2 in \cite{CFCKU}).

The following result from \cite{FMS-1} describes the conditions under which $\overline{X}^{K}_{\theta, b, E}$ is an interpolation space for the couple $(X_{0}, X_{1})$.
\begin{pro} 
Let $\overline{X} = (X_{0}, X_{1})$ be a compatible couple,  $E$  a r.i.  space, $b$ a slowly varying function and $0\leq\theta \leq 1$. Then, the following statements hold
\begin{enumerate}[(a)]
\item  For $t\in \R$ 
\begin{equation}\label{eK}
e^{-\theta t}b(e^{t})K(e^{t},f)
\lesssim \|f\|^{K}_{\theta,b,E}.
\end{equation}
\item  $\overline{X}^{K}_{\theta, b, E}$ is an intermediate space for the
couple $\overline{X}$ if:
\begin{align}
&0< \theta< 1 &\text{ or }\label{c32} \\
& \theta =0 \text{ and } \|b(e^{t})\|_{E(0,\infty)}<\infty&\text{or} \label{c33}\\
& \theta=1 \text{ and } \|b(e^{t})\|_{E(-\infty,0)}<\infty.\label{c34}
\end{align}
\item
 If none of the conditions \eqref{c32} to \eqref{c34} hold, then $\overline{X}^{K}_{\theta, b, E}$ is a trivial space, that is $\overline{X}^{K}_{\theta, b, E}=\{0\}$.
\item
 In any of cases \eqref{c32} to \eqref{c34} the space $\overline{X}^{K}_{\theta, b, E}$ is an
interpolation space for the couple $(X_{0}, X_{1})$.
\end{enumerate}
\end{pro}
Next we introduce the $J$-interpolation method. See \cite{FMS-1}, \cite{FMS-2} and \cite{FMS-3} for more information on $K$ and $J$ spaces.
\begin{dfn}
Let $\overline{A}=(A_0,A_1)$ be a Banach couple, $E$  a r.i.  space, $b$ a slowly varying function and $0\leq\theta \leq 1$. We say that an element $a\in A_0+A_1$ belongs to the space $(A_0,A_1)^J_{\theta,b,E}$ if there exists a representation of $a$ as
\begin{equation}\label{u}
a=\int_{\R} u(e^{t}) dt,
\end{equation}
where $u$ is a strongly measurable function with values in $A_0\cap A_1$, satisfying that
\begin{equation}\label{uJ}
\left\Vert e^{-\theta t} b(e^{t})J(e^{t},u(e^{t}))\right\Vert_{E}<\infty.
\end{equation}
The norm of the element $a$ in the space $(A_0,A_1)^J_{\theta, b,E}$ is given by
\begin{equation*}
\left\Vert a\right\Vert^J_{\theta, b,E}=\inf \Big \{ \left\Vert e^{-\theta t} b(e^{t}) J(e^{t},u(e^{t}))\right\Vert_{E} \Big \},
\end{equation*}
where the infimum is taken over all representations of $a$ satisfying \eqref{u} and \eqref{uJ}.
\end{dfn}
The space $(A_0,A_1)^J_{\theta, b,E}$, $0 < \theta < 1$ is an interpolation space for the couple $(A_{0}, A_{1})$.  
The limit spaces $(A_0,A_1)^J_{0, b,E}$ and $(A_0,A_1)^J_{1, b,E}$ are interpolation spaces for $(A_{0}, A_{1})$ if $\min\{1, e^{-t}\}\tfrac{1}{b(e^{t})}$, $t \in \R$, belongs to $ E'(\R)$, or 
$\min\{1, e^{t}\}\tfrac{1}{b(e^{t})}$, $t \in \R$, is in 
$E'(\R)$ respectively.  
We refer to \cite{B-K} for more details on this question.

In what follows $\overline{A} = (A_{0}, A_{1})$, $\overline{B} = (B_{0}, B_{1})$ and $\overline{C} = (C_{0}, C_{1})$ will be interpolation couples and $T$ will stand for a bounded bilinear operator acting between the couples $\overline{A}$ ,  $\overline{B}$ and $\overline{C}$. More precisely,  there exist constants $k_{i}$, $i=0,1$ such that the operators
\begin{align}
 &T: A_{0} \times B_{0} \longrightarrow C_{0} \\
 &T: A_{1} \times B_{1} \longrightarrow C_{1} \\
\end{align}
satisfy that
$$ \|  T(a,b) \|_{C_{i}} \leq  k_{i} \|   a \|_{A_{i}} \|   b \|_{B_{i}}$$
for all $(a,b) \in A_{i} \times B_{i}$, $i=0,1$.

In these conditions we have the following result.
\begin{thm} \label{Thm JJK}
 Let $0 \leq \theta \leq 1$, $\varphi$ and slowly varying function and $E$ and r.i. space. 
 Then 
  $$T:  \overline{A}_{\theta,m_{\varphi},E}^{J} \times \overline{B}^{K}_{\theta, \varphi, \mathcal{F}(L^{1}, E')} \longrightarrow \overline{C}^{K}_{\theta,\varphi, \mathcal{F}(E,L^{\infty})} $$ 
 is a bounded bilinear operator with norm no greater than $\|   T \| = \max \{\|   T \|_{0}, \|   T_{1} \|\}$.
\end{thm}
\begin{proof}
 Let $(a,b) \in \overline{A}_{\theta,m_{\varphi},E}^{J} \times \overline{B}^{K}_{\theta,\varphi, \mathcal{F}(L^{1}, E')}$. Choose a representation $a = \int_{\R} u(e^{t}) dt $ 
 satisfying that 
 \begin{equation*}
 \|   e^{-\theta t} \varphi(e^{t}) J(e^{t}, u(e^{t})) \|_{E} \leq (1+ \varepsilon) \|  a  \|^{J}_{\theta, \varphi, E} .
\end{equation*}
Similarly, for each $t>0$ we choose decompositions of $b \in B_{0} + B_{1}$ satisfying that
\begin{equation}
 \|   b_{0}(t) \|_{B_{0}} + t \|   b_{1}(t) \|_{B_{1}} \leq (1 + \varepsilon) K(t,b).
\end{equation}
Since
\begin{equation*}
 T(a,b) = T \Big(  \int_{\R} u(e^{t}) dt, b   \Big)=  \int_{\R} T \big (u(e^{t}) , b   \big) dt,
\end{equation*}
it is not difficult to check that
\begin{align*}
 \frac{\varphi(e^{y})}{e^{\theta y}} &K \big ( e^{y}, T(a,b) \big ) 
 \leq \frac{\varphi(e^{y})}{e^{\theta y}} \int_{\R} K \big (e^{y}, T \big (u(e^{x}), b \big) \big ) dx \\
 & \leq \|   T \| \int_{\R} \frac{\varphi(e^{y})}{e^{\theta y}} \Big [  \| u(e^{x})   \|_{0}  \, \|   b_{0}(e^{y-x}) \|_{0} 
 + e^{y}   \| u(e^{x})   \|_{1} \, \|   b_{1}(e^{y-x}) \|_{1}      \Big ] dx \\
 &\leq \|   T \| \int_{\R} \frac{\varphi(e^{y})}{e^{\theta y}} J \big (e^{x}, u(e^{x}) \big) \, (1+ \varepsilon) K \big( e^{y-x},b \big ) dx \\
 &\leq (1+ \varepsilon) \|   T \| \int_{\R} \frac{m_{\varphi}(e^{x})}{e^{\theta x}} J(e^{x}, u(e^{x})) \frac{\varphi(e^{y-x})}{e^{\theta(y-x)}} K(e^{y-x},b) dx \\
 & = (1+ \varepsilon) \|   T \| \Big (  \frac{m_{\varphi}(e^{x})}{e^{\theta x}} J(e^{x}, u(e^{x}))   \Big) 
 \ast \Big (  \frac{\varphi(e^{x})}{e^{\theta(x)}} K(e^{x},b) \Big) (y).
\end{align*}
Hence, taking infimum in $\varepsilon >0$ 
\begin{equation*}
 \Big \|  \frac{\varphi(e^{y})}{e^{\theta y}} K(e^{y}, T(a,b))  \Big \|_{E} \leq  \|   T \|  \Big \|  \Big (  \frac{m_{\varphi}(e^{x})}{e^{\theta x}} J(e^{x}, u(e^{x}))   \Big) 
 \ast \Big (  \frac{\varphi(e^{x})}{e^{\theta(x)}} K(e^{x},b) \Big) (y) \Big \|_{E} .
\end{equation*}
Now, apply Thm.~\ref{GYI-thm} to obtain that
\begin{equation*}
 \|   T(a,b) \| \leq \|   T \| \, \|   a \|^{J}_{\theta, m_{\varphi}, E} \, \|   b \|^{K}_{\theta, \varphi, \mathcal{F}(L^{1},E')}.
\end{equation*}
This completes the proof.
\end{proof}
Also, we can establish a bilinear interpolation theorem for $J$-spaces.
\begin{thm}
 Let $0 \leq \theta \leq 1$, $\varphi$ and slowly varying function and $E$ and r.i. space. 
 Then, for any interpolation functor $\mathcal{F}$, 
 $$T:  \overline{A}_{\theta,\varphi,E}^{J} \times \overline{B}^{J}_{\theta, m_{\varphi}, \mathcal{F}(L^{1}, E')} \longrightarrow \overline{C}^{J}_{\theta,\varphi, \mathcal{F}(E,L^{\infty})} $$
 is a bounded bilinear operator with norm no greater than $\|   T \| = \max \{\|   T \|_{0}, \|   T_{1} \|\}$.
\end{thm}
\begin{proof}
 Let $(a,b) \in  \overline{A}_{\theta,\varphi,E}^{J} \times \overline{B}^{J}_{\theta, m_{\varphi}, \mathcal{F}(L^{1}, E')}$ and choose representations of both elements
 $$a = \int_{\R} u(e^{x}) dx \quad \text{ and }  \quad b = \int_{\R} v(e^{y}) dy $$
satisfying that 
\begin{align*}
\Big \|   e^{-\theta x} \varphi(e^{(x)}) J(e^{(x)}, u(e^{(x)})) \Big \|_{E}   &\leq (1+ \varepsilon)  \|   a \|_{\theta,\varphi,E}^{J} \\
\Big \|   e^{-\theta y} \varphi(e^{(y)}) J(e^{(y)}, u(e^{(y)})) \Big \|_{E}   &\leq (1+ \varepsilon)  \|   b \|_{\theta,m_{\varphi}, \mathcal{F}(L^{1}, E') }^{J}.
\end{align*}
Since 
$$T(a,b) = T \Big( \int_{\R} u(e^{x}) dx,    \int_{\R} v(e^{y}) dy  \Big) = \int_{\R} \int_{\R } T \big (  u(e^{x}), v(e^{y-x}) \big) dx dy,$$ we claim that  
$w(e^{y}) = \int_{\R} T \big (  u(e^{x}), v(e^{y-x}) \big) dx \in C_{0} \cap C_{1}$
and that we may consider $\int_{\R} w(e^{y})dy$  a representation of $T(a,b)$ as an element of 
$\overline{C}^{J}_{\theta,\varphi, \mathcal{F}(E,L^{\infty})}$.
In fact, the $J$-functional of $w(e^{y})$, $y \in \R$, satisfies the inequalities
\begin{align*}
 J(e^{y}, w(e^{y})) &= J \Big  (   e^{y}, \int_{\R}   T \big (  u(e^{x}), v(e^{y-x}) \big)  dx  \Big ) \\
 & \leq \int_{\R}  J \Big(  e^{y}, T \big (  u(e^{x}), v(e^{y-x}) \big) \Big) dx \\
 & \leq \max  \big \{ \|   T \|_{0}, \|   T \|_{1}  \big \}   \int_{\R}     J \big ( e^{x},u(e^{x}) \big )  J \big ( e^{y-x},u(e^{y-x}) \big )     dx
\end{align*}
where the last integral can be estimated by a convolution as follows
\begin{equation} \label{J1}
 \begin{aligned}
 \int_{\R}  &   J \big ( e^{x},u(e^{x}) \big )  \; J \big ( e^{y-x},v(e^{y-x}) \big )     dx  \\
& \leq    \int_{\R}    e^{-\theta x}  \varphi(e^{x}) J \big ( e^{x},u(e^{x}) \big )  \frac{e^{\theta y}}{\varphi(e^{y})} e^{-\theta(y-x))}m_{\varphi}(e^{y-x}))  J \big ( e^{y-x},v(e^{y-x}) \big ) dx  \\
&= \frac{e^{\theta y}}{\varphi(e^{y})}   \Bigg [ \bigg(e^{-\theta x}  \varphi(e^{x}) J \big ( e^{x},u(e^{x}) \big ) \bigg) \ast  
 \bigg( e^{-\theta x}m_{\varphi}(e^{x})  J \big ( e^{x},v(e^{x}) \big ) \bigg)    \Bigg ](y).
\end{aligned}
\end{equation}
This puts us in a position to apply  \eqref{GYI} to obtain that
\begin{align*}
 &\|  e^{- \theta y} \varphi(e^{y})   J(e^{y}, w(e^{y}))  \|_{\mathcal{F}(E, L^{\infty})} \\ 
 & \leq
 \bigg \|   \bigg [ \Big(e^{-\theta x}  \varphi(e^{x}) J \big ( e^{x},u(e^{x}) \big ) \Big) \ast  
 \Big (e^{-\theta x}m_{\varphi}(e^{x})  J \big ( e^{x},v(e^{x}) \big ) \Big)    \bigg ] \bigg\|_{\mathcal{F}(E, L^{\infty})} \\
 & \leq \Big\|   e^{-\theta x}m_{\varphi}(e^{x})  J \big ( e^{x},u(e^{x}) \big )    \Big\|_{E}     \;
\Big \|  e^{-\theta x}m_{\varphi}(e^{x})  J \big ( e^{x},v(e^{x}) \big )   \Big \|_{\mathcal{F}(L^{1}, E^{'}))} \\
& \leq (1 + \varepsilon)^{2} \|   a \|_{\theta, \varphi, E}  \; \|   b \|_{\theta, m_{\varphi}, \mathcal{F}(L^{1}, E^{'}))}.
\end{align*}
Therefore, the function $y \leadsto e^{-\theta y} \varphi(e^{y}) J(e^{y},w(e^{y}))$, $y \in \R$, belongs to $F(E, L^{\infty})$ and so $J(e^{y},w(e^{y}))$, $y \in \R$, is finite \emph{a.e.} in $\R$. This not only proves that 
$$T(a,b) =\int_{\R} w(e^{y})dy$$
is a representation of $T(a,b)$ as an element of $\overline{C}^{J}_{\theta,\varphi, \mathcal{F}(E,L^{\infty})}$  but also establishes that 
$$\| T(a,b) \|_{\overline{C}^{J}_{\theta,\varphi, \mathcal{F}(E,L^{\infty})}} \leq \max \big \{\|   T \|_{0}, \|   T_{1} \| \big \} \; \| a \|_{\overline{A}_{\theta,\varphi,E}^{J}} \| b \|_{\overline{B}^{J}_{\theta, m_{\varphi}, \mathcal{F}(L^{1}, E')}}.$$
This completes the proof.
\end{proof}

\section{Aplications to bilinear multipliers} \label{bilinear multipliers}
We ilustrate the usefulness of our results by applying them in the context of bilinear multipliers. More precisely, we work on a result by O. Blasco on spaces of bilinear multipliers in \cite{Blasco2009}.

Let $f$ and $g$ be periodic functions defined on the torus $\mathbb{T}$. The operator $\mathcal{P}_{m}$ is defined as 
\begin{equation*}
 \mathcal{P}_{m} (f,g) (\theta) = \sum_{k, k' \in \Z} \hat{f}(k) \hat{g}(k') \, m_{k,k'} \,  e^{2 \pi i \theta (k + k')} .
\end{equation*}
An easy duality argument shows that
$$\mathcal{P}_{m} : L^{p_{1}}(\mathbb{T}) \times  L^{p_{2}}(\mathbb{T}) \longrightarrow L^{p_{3}}(\mathbb{T})$$
if and only if there exists $C>0$ such that
$$ \Big | \sum_{k, k' \in \Z} \hat{f}(k) \hat{g}(k') \, m_{k,k'} \, \hat{h}(k + k')    \Big |  \leq C \|   f \|_{p_{1}}   \|  g  \|_{p_{2}}  \|   h \|_{p_{3}'}$$
for all $h \in L^{p_{3}'}(\mathbb{T})$. In this case we say that $m = (m_{k,k'})_{\Z^{2}}$ belongs to $\mathcal{BM}_{(p_{1}, p_{2}, p_{3})}$, the space of bilinear multipliers of type $(p_{1}, p_{2}, p_{3})$.

 Following Thm.~2.8 of \cite{Blasco2009}, 
it is easy to establish that for $m=(m_{k,k'})_{\Z^{2}} \in \ell{p}(\Z^{2})$, the operators
\begin{align*}
 \mathcal{P}_{m} &: L^{p}(\mathbb{T}) \times  L^{p}(\mathbb{T}) \longrightarrow L^{\infty}(\mathbb{T}) \\
\mathcal{P}_{m} &: L^{p}(\mathbb{T}) \times  L^{1}(\mathbb{T}) \longrightarrow L^{p'}(\mathbb{T})
\end{align*}
are bounded.  In this position, we can apply Thm. \ref{Thm JJK} to obtain that
$$\mathcal{P}_{m} : (L^{p}, L^{p})_{0, m_{\ell^{-1/p}(t)},L^{1}}^{J} \times 
 (L^{p}, L^{1})_{0, \ell^{-1/p}(t),L^{\infty}}^{K} \longrightarrow (L^{\infty}, L^{p'})_{0, \ell^{-1/p}(t),L^{\infty}}^{K}.$$
 Now we proceed to identify these spaces. 
 
 The space $(L^{p}, L^{1})_{0, \ell^{-1/p}(t),L^{\infty}}$ is the grand-Lebesgue space $L^{p)}$  introduced by Iwaniec and Sbordone in   \cite{Iwaniec1992}  in connection with the
study of the integrability properties of the jacobian determinant. It is defined by the norm
 $$\|f\|_{L^{p)}}=\|\ell^{-1/p}(t)\|f^{**}(s)\|_{L_p(t,2 \pi)}\|_{L_\infty(0,2 \pi)}<\infty.$$
See  \cite{Fiorenza2004} for additionally information of the norm of $L^{p)}(\Omega)$.
On the other hand,  the space $(L^{\infty}, L^{p'})_{0, \ell^{-1/p}(t),L^{\infty}}$ is continuously embedded in $L_{exp}$. Indeed, 
\begin{align*}
(L^{\infty}, L^{p'})_{0, \ell^{-1/p}(t), L^{\infty}} &= ( L^{p'} , L^{\infty})_{1, \ell^{-1/p}(t), L^{\infty}} \\
& = \Big( (L^{1}, L^{\infty})_{1 - 1/p', 1, L^{p'} , L^{\infty}}  \Big)_{ 1, \ell^{-1/p}(t), L^{\infty} } \\
& = (L^{1}, L^{\infty})_{1, \ell^{1/p}(t^{1/p'}), L^{\infty}} \hookrightarrow L_{exp}.
\end{align*}
This shows that 
$$\mathcal{P}_{m} : L^{p}\times L^{p)} \longrightarrow L_{exp} , $$
which improves the domain of the original operator out of the $L^{p} scale$ and restrict the range.
\bibliographystyle{amsplain}
\bibliography{/Users/Pedro/Dropbox/Tex/bibliography} 
\linespread{1}

\end{document}